\documentclass[11pt]{article}
\usepackage[a4paper,margin=29mm]{geometry}
\usepackage{amsmath,amssymb,amsthm,mathtools}
\usepackage{microtype}
\usepackage[hidelinks]{hyperref}
\usepackage{enumitem}
\usepackage{tikz}
\usepackage{float}
\usetikzlibrary{angles,quotes}

\newtheorem{theorem}{Theorem}[section]
\newtheorem{lemma}[theorem]{Lemma}
\newtheorem{proposition}[theorem]{Proposition}
\newtheorem{corollary}[theorem]{Corollary}
\theoremstyle{remark}
\newtheorem{remark}[theorem]{Remark}
\newtheorem{definition}[theorem]{Definition}

\DeclareMathOperator{\atanTwo}{atan2}

\newcommand{\areg}{\alpha_*}
\newcommand{\Kclass}{\mathcal K}
\newcommand{\Oclass}{\mathcal O}
\newcommand{\Iclass}{\mathcal I}

\title{On Exact Dihedral Angle Sum Ranges for\\Orthocentric and Isogonic Tetrahedra}
\author{Sergey Korotov\\
\small Department of Business and Mathematics, M\"alardalen University\\
\small Box 883, 721 23 V\"aster\aa s, Sweden\\
\small \texttt{sergey.korotov@mdu.se}}
\date{19 October 2026}

\begin{document}
\maketitle
\vspace{-0.8em}

\begin{abstract}
For a nondegenerate Euclidean tetrahedron $T$, let $\Sigma(T)$ denote the sum of its six interior dihedral angles. We establish the exact range
\[
 \left[6\arccos\frac13,\,3\pi\right)
\]
for the orthocentric and isogonic classes. The orthocentric part gives a short alternative proof of the sharp lower bound; the main new result is the corresponding isogonic theorem. Starting from classical insphere-contact geometry, we derive a reciprocal half-angle identity and a sharp defect inequality controlled by the scale-invariant dispersion of the reciprocal tangent-length parameters. This quantitative refinement shows that any nonzero dispersion forces a positive dihedral-angle-sum gap, and that near equality forces the normalized reciprocal parameters to be close to equality. Finally, the common kite subfamily realizes every value in the interval and proves sharpness of both endpoints.
\end{abstract}



\section{Introduction}
Let $T$ be a nondegenerate Euclidean tetrahedron and write
\[
 \Sigma(T)=\sum_{e\subset T}\theta_e
\]
for the sum of its six interior dihedral angles. The classical sharp universal bounds are
\begin{equation}
 2\pi<\Sigma(T)<3\pi,
 \label{eq:universal}
\end{equation}
see Jarden~\cite{Jarden1941}, Gaddum~\cite{Gaddum1952}, and the geometric treatment in Sharygin~\cite{Sharygin1986}. The regular tetrahedron has
\begin{equation}
 \Sigma_{\rm reg}=6\arccos\frac13.
 \label{eq:regular-value}
\end{equation}

We denote by $\Oclass$, $\Iclass$, and $\Kclass$ the orthocentric, isogonic, and tetrahedral-kite classes, respectively. A tetrahedron is \emph{orthocentric} if its four altitudes are concurrent. Brown~\cite{Brown1924} uses the older term \emph{orthogonal} for this class. Classical equivalent characterizations include pairwise perpendicular opposite edges,
\[
 AB^2+CD^2=AC^2+BD^2=AD^2+BC^2,
\]
and the additive squared-edge representation
\begin{equation}
 d_{ij}^{\,2}=a_i+a_j,
 \qquad i\ne j,
 \label{eq:orthometric}
\end{equation}
with real parameters for which all pair sums are positive; see Edmonds--Hajja--Martini~\cite{EdmondsHajjaMartini2005} and Sharygin~\cite[Problem~312]{Sharygin1986}.

A tetrahedron is \emph{isogonic} if the four cevians joining its vertices to the contact points of the insphere with the opposite faces are concurrent. Older sources use the term \emph{isogonal}~\cite{Brown1924}, while Hajja--Walker use \emph{inspherical}~\cite{HajjaWalker2004}. Hajja--Hammoudeh~\cite[Theorem~4]{HajjaHammoudeh2014} give the equivalent metric characterization
\begin{equation}
 d_{ij}^{\,2}=p_i^2+p_ip_j+p_j^2,
 \qquad p_i>0.
 \label{eq:isogonicmetric}
\end{equation}
Every positive parameter quadruple is admissible; a positive volume formula is recalled in Section~5. Multiplying all $p_i$ by the same positive constant produces a similar tetrahedron, so the parameter quadruple represents a similarity class modulo common scaling.

A \emph{tetrahedral kite}, also called a \emph{kite simplex} in the general simplex literature, is a tetrahedron with an equilateral face and with the fourth vertex equidistant from the three vertices of that face. Equivalently, after relabelling,
\[
 d_{12}=d_{13}=d_{23},
 \qquad
 d_{14}=d_{24}=d_{34}.
\]
Within the isogonic parametrization \eqref{eq:isogonicmetric}, this is exactly the case in which three of the four parameters are equal. An arbitrary equilateral-base tetrahedron need not be a kite, because its three lateral edges need not be equal.

Further classical and modern treatments of the special families above include Mandan~\cite{Mandan1961}, Hajja--Hayajneh--Martini~\cite{HajjaHayajnehMartini2016}, and the Russian sources Skopets--Ponarin~\cite{SkopetsPonarin1974}, Sharygin~\cite{Sharygin1986}, and Ponarin~\cite{Ponarin2009}.

Hammoudeh~\cite{Hammoudeh2026} recently obtained the exact orthocentric range and, in Appendix~A of the cited Version~2.0.0, recorded the corresponding regular lower bound for isogonic tetrahedra as unresolved; Conjecture~A.2 proposes a stronger ordered reduction to kites. The new contributions of the present paper are the isogonic lower bound and exact range, an independent geometric derivation of the reciprocal half-angle identity, and a quantitative defect estimate. The orthocentric part is an alternative proof, while the kite-intersection theorem is recalled and reproved in dimension three because its parametrization supplies sharpness.

\begin{theorem}[Exact isogonic range]\label{thm:isogonic-main}
The set of dihedral angle sums of nondegenerate isogonic tetrahedra is
\[
 \{\Sigma(T):T\ \text{is isogonic}\}
 =\left[6\arccos\frac13,\,3\pi\right).
\]
The lower endpoint is attained only by the regular tetrahedron.
\end{theorem}

For comparison, we also give a short proof of the sharp orthocentric lower bound. Combined with the universal upper bound and the kite calculation below, this recovers the exact orthocentric range without using the cited manuscript as a logical input. The orthocentric proof is based on a four-corner conservation law and a supporting-line estimate. The isogonic proof has a different mechanism: classical $120^\circ$ insphere-contact geometry yields a reciprocal half-angle identity, and two exact symmetry-defect identities lead to the quantitative estimate.

The two classes are linked by kites. Hajja--Hammoudeh--Hayajneh~\cite{HajjaHammoudehHayajneh2015} proved, more generally, that the intersection of any two among the orthocentric, circumscriptible, isodynamic, and isogonic simplex families is precisely the kite family. In dimension three the circumscriptible class is also known as the 4-ball or balloon family; its dihedral angle sum has the same exact interval in Korotov--K\v{r}\'i\v{z}ek~\cite[Theorem~9]{KorotovKrizek2026}. We include a short direct verification of
\begin{equation}
 \Oclass\cap\Iclass=\Kclass
 \label{eq:intersection-intro}
\end{equation}
(up to relabelling), because the same parametrization supplies the sharpness construction. The common kite subfamily does not prove either lower bound, but once the two lower bounds are known it realizes every value up to the universal upper endpoint.

Sharp or exact class-dependent bounds are known for several other tetrahedral families. For equifacial (isosceles) tetrahedra, Katsuura~\cite[Theorem~3]{Katsuura2020} proved
\[
 2\pi<\Sigma(T)\le 6\arccos\frac13,
\]
with equality on the right only for the regular tetrahedron; thus the regular tetrahedron is a maximizer in that class, in contrast with the orthocentric and isogonic classes considered here. The same paper also treats right (trirectangular) tetrahedra. The exact path-tetrahedron range is $(2\pi,5\pi/2)$~\cite[Theorem~2]{KorotovKrizek2026}, while Brandts--Cihangir--K\v{r}\'i\v{z}ek~\cite{BrandtsCihangirKrizek2015} established tight angle-sum bounds for nonobtuse simplices. The purity of orthocentric corner angles is discussed by Hajja--Hayajneh~\cite{HajjaHayajneh2014}.

Dihedral angles and solid-angle-based quantities are standard local indicators of tetrahedral shape quality in mesh generation and finite element computations; see, for example, Liu--Joe~\cite{LiuJoe1994}. Exact restrictions on the aggregate quantity $\Sigma(T)$ over geometrically distinguished tetrahedral families provide structural information complementary to the usual minimum- and maximum-angle conditions. The results below are geometric rather than algorithmic, but they delimit the angle-sum behavior that is intrinsically possible within two classical tetrahedral classes.

\section{Corner measures and the dihedral angle sum}
Consider a trihedral corner whose three unit edge vectors have pairwise scalar products $u,v,w$. The square root of their Gram determinant,
\begin{equation}
 P=\sqrt{1-u^2-v^2-w^2+2uvw}>0,
 \label{eq:polar-sine}
\end{equation}
is commonly called the \emph{polar sine}. Its solid-angle measure is given by Eriksson's formula~\cite{Eriksson1990},
\begin{equation}
 \Omega=2\,\atanTwo(P,1+u+v+w).
 \label{eq:solid-angle}
\end{equation}
Here $\atanTwo(y,x)\in(0,\pi)$ denotes the polar angle of $(x,y)$ for $y>0$; in particular, for fixed $y>0$ it is strictly decreasing as a function of $x\in\mathbb R$.
For the pure corners needed below, define the polar measure $\beta$ by
\begin{equation}
 \beta=
 \begin{cases}
  \arcsin P, & u,v,w>0,\\[2mm]
  \pi/2, & u=v=w=0,\\[2mm]
  \pi-\arcsin P, & u,v,w<0.
 \end{cases}
 \label{eq:polar-measure}
\end{equation}
Thus $\sin\beta=P$ in all three cases.

Let $\Omega_i$ be the solid-angle measure at vertex $A_i$. The spherical link at $A_i$ has area
\[
 \Omega_i=\sum_{e\ni A_i}\theta_e-\pi.
\]
Summing over the four vertices counts each dihedral angle twice and gives
\begin{equation}
 \Sigma(T)=2\pi+\frac12\sum_{i=1}^4\Omega_i.
 \label{eq:solid-sum}
\end{equation}
We shall also use the following established orthocentric identity of Hajja--Hammoudeh~\cite[Theorem~1(b)]{HajjaHammoudeh2014}:
\begin{equation}
 \beta_1+\beta_2+\beta_3+\beta_4=\pi
 \qquad\text{for every nondegenerate orthocentric tetrahedron.}
 \label{eq:polar-sum}
\end{equation}
The purity property is characteristic of the orthocentric setting among the classical special families; see Hajja--Hayajneh~\cite{HajjaHayajneh2014}. Indeed, from \eqref{eq:orthometric} the cosine of a planar angle at vertex $i$ has the form
\[
 \frac{a_i}{\sqrt{(a_i+a_j)(a_i+a_k)}},
\]
so all three planar angles at a fixed vertex have cosines of the same sign. Hence every corner of an orthocentric tetrahedron is pure.

\section{A short alternative proof for orthocentric tetrahedra}
Put
\begin{equation}
 \areg:=\arccos\frac13,
 \qquad
 c_*:=3\areg-\pi,
 \qquad
 m_*:=\frac{2\sqrt2}{3},
 \label{eq:constants}
\end{equation}
and define the affine function
\begin{equation}
 L(\beta):=c_*+m_*\left(\beta-\frac{\pi}{4}\right).
 \label{eq:support-line}
\end{equation}

\begin{lemma}[Supporting-line corner estimate]\label{lem:corner-affine}
Every pure nondegenerate trihedral corner satisfies
\begin{equation}
 \Omega\ge L(\beta).
 \label{eq:corner-affine}
\end{equation}
Equality holds if and only if the corner is regular, i.e. its three planar angles are all $\pi/3$; in that case $\beta=\pi/4$.
\end{lemma}

\begin{proof}
First suppose the corner is acute, so $u,v,w>0$ and $0<\beta<\pi/2$. Put $s=u+v+w$ and
\[
 \phi(q):=(1-q)^2(1+2q)=1-3q^2+2q^3,
 \qquad 0\le q\le1.
\]
The function $\phi$ decreases strictly from $1$ to $0$ on $[0,1]$. Hence for each $\beta\in(0,\pi/2]$ there is a unique $q=q(\beta)\in[0,1)$ such that
\begin{equation}
 \sin^2\beta=\phi(q).
 \label{eq:qbeta}
\end{equation}
Here $q(\pi/2)=0$, while $q(\beta)\uparrow1$ as $\beta\downarrow0$.

The elementary inequalities
\[
 u^2+v^2+w^2\ge\frac{s^2}{3},
 \qquad
 uvw\le\left(\frac{s}{3}\right)^3
\]
give
\[
 \sin^2\beta
 =1-u^2-v^2-w^2+2uvw
 \le 1-\frac{s^2}{3}+\frac{2s^3}{27}
 =\phi\!\left(\frac{s}{3}\right).
\]
Since $\phi$ is decreasing, $s\le3q$. Using $P=\sin\beta$ in \eqref{eq:solid-angle} and the monotonicity of $\atanTwo(\sin\beta,X)$ in $X$,
\begin{align}
 \Omega
 &=2\atanTwo(\sin\beta,1+s)\notag\\
 &\ge2\atanTwo(\sin\beta,1+3q)
 =:f(\beta).                                           \label{eq:f-atan}
\end{align}
The right-hand side is the solid-angle content of the equiangular corner with common planar cosine $q$. By the spherical cosine law, each of its three dihedral angles equals $\arccos(q/(1+q))$, and therefore
\begin{equation}
 f(\beta)
 =3\arccos\!\left(\frac{q(\beta)}{1+q(\beta)}\right)-\pi.
 \label{eq:f-beta}
\end{equation}
Equality in \eqref{eq:f-atan}--\eqref{eq:f-beta} occurs exactly when $u=v=w$.

Differentiating the parametrization \eqref{eq:qbeta} gives
\begin{equation}
 f'(\beta)=\frac{\sqrt{3-2q}}{1+q},
 \qquad
 f''(\beta)=\frac{(4-q)\sqrt{1+2q}}{3(1+q)^2}>0.
 \label{eq:f-derivatives}
\end{equation}
Thus $f$ is strictly convex. At $\beta=\pi/4$ one has $q=1/2$, and hence
\[
 f\!\left(\frac{\pi}{4}\right)=c_* ,
 \qquad
 f'\!\left(\frac{\pi}{4}\right)=m_*.
\]
The tangent-line inequality for a convex function gives $f(\beta)\ge L(\beta)$, proving \eqref{eq:corner-affine} for acute corners. Equality requires $u=v=w$ and $\beta=\pi/4$, hence $u=v=w=1/2$.

For the orthogonal corner $u=v=w=0$, one has $\Omega=\beta=\pi/2$. For an obtuse pure corner write
\[
 u=-U,\qquad v=-V,\qquad w=-W,
 \qquad 0<U,V,W<1,
\]
and put $\rho:=U+V+W$. Since $UVW\le UV+UW+VW$,
\[
 1-P^2
 =U^2+V^2+W^2+2UVW
 \le U^2+V^2+W^2+2UV+2UW+2VW
 =\rho^2.
\]
Because $\beta=\pi-\arcsin P$,
\[
 \rho\ge\sqrt{1-P^2}=-\cos\beta.
\]
Using \eqref{eq:solid-angle} and the monotonicity of $\atanTwo(P,X)$ in $X$,
\begin{equation}
 \Omega
 =2\atanTwo(P,1-\rho)
 \ge2\atanTwo(\sin\beta,1+\cos\beta)
 =\beta.
 \label{eq:obtuse-beta}
\end{equation}
Strict convexity above gives $L(\pi/2)<f(\pi/2)=\pi/2$, and $L'(\beta)=m_*<1$. Therefore $L(\beta)<\beta$ for every $\beta\ge\pi/2$. Hence \eqref{eq:corner-affine} holds strictly for orthogonal and obtuse corners.
\end{proof}

\begin{theorem}[Sharp orthocentric lower bound]\label{thm:orth-lower}
Every nondegenerate orthocentric tetrahedron satisfies
\begin{equation}
 \Sigma(T)\ge6\arccos\frac13.
 \label{eq:orth-lower}
\end{equation}
Equality holds if and only if $T$ is regular.
\end{theorem}

\begin{proof}
Apply Lemma~\ref{lem:corner-affine} at the four vertices and use \eqref{eq:polar-sum}:
\[
 \sum_{i=1}^4\Omega_i
 \ge4c_*+m_*\left(\sum_{i=1}^4\beta_i-\pi\right)
 =4c_*.
\]
Now \eqref{eq:solid-sum} gives
\[
 \Sigma(T)
 \ge2\pi+2c_*
 =6\arccos\frac13.
\]
If equality holds, equality must hold in the corner estimate at every vertex. Hence every planar face angle is $\pi/3$, so all faces are equilateral and $T$ is regular. The converse is immediate.
\end{proof}

\begin{remark}
The same sharp lower bound and exact orthocentric range were independently obtained by a different argument in \cite{Hammoudeh2026}. The proof above uses the earlier constant polar-measure identity \eqref{eq:polar-sum} of Hajja--Hammoudeh~\cite{HajjaHammoudeh2014} as its only external orthocentric input; the full range follows below from the kite family and the classical universal upper bound.
\end{remark}

\section{Kites: the common bridge}
\begin{definition}
A \emph{tetrahedral kite} (hereafter simply a \emph{kite}) is a tetrahedron with an equilateral face and with the fourth vertex equidistant from the three vertices of that face. Equivalently, after relabelling, its three base edges are equal and its three lateral edges are equal. We denote this family by $\Kclass$.
\end{definition}

Normalize the base edge to $1$ and let $t$ be the common lateral edge. Nondegeneracy is equivalent to
\begin{equation}
 t>\frac1{\sqrt3}.
 \label{eq:kite-domain}
\end{equation}

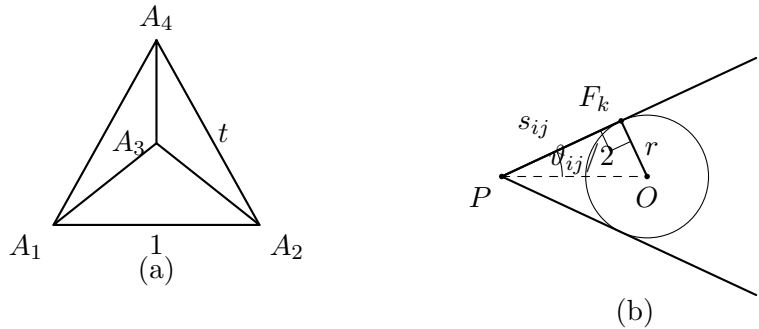
\begin{figure}[H]
\centering
\begin{tikzpicture}[line cap=round,line join=round,scale=0.80]
\begin{scope}[xshift=-4.5cm]
 \coordinate (A1) at (-1.7,-0.8);
 \coordinate (A2) at (1.7,-0.8);
 \coordinate (A3) at (0,0.55);
 \coordinate (A4) at (0,2.25);
 \draw[thick] (A1)--(A2)--(A3)--cycle;
 \draw[thick] (A4)--(A1) (A4)--(A2) (A4)--(A3);
 \node[below left] at (A1) {$A_1$};
 \node[below right] at (A2) {$A_2$};
 \node[left] at (A3) {$A_3$};
 \node[above] at (A4) {$A_4$};
 \path (A1)--(A2) node[midway,below] {$1$};
 \path (A4)--(A2) node[midway,right] {$t$};
 \node at (0,-1.55) {(a)};
\end{scope}
\begin{scope}[xshift=1.2cm]
 \coordinate (P) at (0,0);
 \coordinate (U) at (4.2,1.958);
 \coordinate (L) at (4.2,-1.958);
 \coordinate (O) at (2.4,0);
 \coordinate (F) at (1.971,0.919);
 \draw[thick] (P)--(U) (P)--(L);
 \draw[dashed] (P)--(O);
 \draw (O) circle[radius=1.014];
 \draw[thick] (P)--(F) node[midway,above left] {$s_{ij}$};
 \draw[thick] (O)--(F) node[midway,right] {$r$};
 \fill (P) circle (1.3pt) node[below left] {$P$};
 \fill (O) circle (1.3pt) node[below] {$O$};
 \fill (F) circle (1.3pt) node[above left] {$F_k$};
 \pic [draw, "$\vartheta_{ij}/2$", angle eccentricity=1.35, angle radius=8mm] {angle=O--P--F};
 \pic [draw, angle radius=3mm] {right angle=P--F--O};
 \node at (2.2,-2.25) {(b)};
\end{scope}
\end{tikzpicture}
\caption{(a) A normalized tetrahedral kite, with base edge $1$ and common lateral edge $t$. (b) The perpendicular cross-section used in Lemma~\ref{lem:geometric-half}: the edge $A_iA_j$ meets the section at $P$, the incenter is $O$, and $F_k$ is a tangency point.}
\label{fig:kite-crosssection}
\end{figure}

\begin{proposition}[Kite intersection; cf. \cite{HajjaHammoudehHayajneh2015}]\label{prop:intersection}
Up to relabelling,
\begin{equation}
 \Oclass\cap\Iclass=\Kclass.
 \label{eq:intersection}
\end{equation}
\end{proposition}

\begin{proof}
Although \eqref{eq:intersection} follows from the general theorem in \cite{HajjaHammoudehHayajneh2015}, we include the elementary three-dimensional proof because the same parametrization also yields the complete kite range needed for sharpness.

Every kite is orthocentric: in \eqref{eq:orthometric} take the three base parameters equal to $1/2$ and the apex parameter equal to $t^2-1/2$.

Every kite is also isogonic. In \eqref{eq:isogonicmetric} take three parameters equal to $b>0$ and the apex parameter equal to $sb$, $s>0$. The squared base edge is $3b^2$ and the squared lateral edge is $b^2(s^2+s+1)$, hence
\begin{equation}
 t^2=\frac{s^2+s+1}{3}.
 \label{eq:kite-s}
\end{equation}
As $s$ ranges over $(0,\infty)$, the right-hand side ranges over $(1/3,\infty)$, exactly the nondegenerate kite range \eqref{eq:kite-domain}.

Conversely, suppose an isogonic tetrahedron \eqref{eq:isogonicmetric} is orthocentric. Equality of the three sums of squared opposite edges gives
\begin{align}
 (p_1-p_4)(p_2-p_3)&=0,\label{eq:prod1}\\
 (p_1-p_3)(p_2-p_4)&=0,\label{eq:prod2}\\
 (p_1-p_2)(p_3-p_4)&=0.\label{eq:prod3}
\end{align}
The set of these three relations is invariant under relabelling. If no three parameters were equal, their equality partition would have type $1+1+1+1$, $2+1+1$, or $2+2$. The first type contradicts any one of \eqref{eq:prod1}--\eqref{eq:prod3}. In the second type, relabel so that $p_1=p_2=a$ while $p_3=b$ and $p_4=c$ are distinct from $a$ and from one another; then \eqref{eq:prod1} gives $(a-c)(a-b)\ne0$, a contradiction. In the third type, relabel so that $p_1=p_2=a$ and $p_3=p_4=b$ with $a\ne b$; then \eqref{eq:prod2} gives $(a-b)^2\ne0$. Hence at least three parameters are equal.

If, say, $p_1=p_2=p_3=b$ and $p_4=a$, then
\[
 d_{12}=d_{13}=d_{23}=\sqrt3\,b,
 \qquad
 d_{14}=d_{24}=d_{34}=\sqrt{a^2+ab+b^2}.
\]
Thus the face $A_1A_2A_3$ is equilateral and the three lateral edges are equal, so the tetrahedron is a kite.
\end{proof}

\begin{remark}[Neither class contains the other]
Neither inclusion holds. For
\[
 A=(0,0,0),\quad B=(1,0,0),\quad C=(0,2,0),\quad D=(0,0,3),
\]
the squared edge lengths are $1,4,9,5,10,13$ and the three opposite-pair sums are all $14$. Hence the tetrahedron is orthocentric, but its six edge lengths are distinct, so it is not a kite and therefore not isogonic. Conversely, the isogonic parameters $(p_1,p_2,p_3,p_4)=(1,1,2,3)$ give squared edges $3,7,13,7,13,19$, whose opposite-pair sums are $22,20,20$; the resulting tetrahedron is not orthocentric. Thus
\[
 \Oclass\not\subseteq\Iclass,
 \qquad
 \Iclass\not\subseteq\Oclass.
\]
\end{remark}

For later use we record the complete kite calculation. Direct scalar products of face normals give
\begin{equation}
 K(t):=\Sigma(T_t)
 =3\left[
 \arccos\frac{1}{\sqrt{3(4t^2-1)}}
 +\arccos\frac{2t^2-1}{4t^2-1}
 \right].
 \label{eq:kite-sum}
\end{equation}
Differentiation simplifies to
\begin{equation}
 K'(t)=\frac{6(t-1)}{(4t^2-1)\sqrt{3t^2-1}}.
 \label{eq:kite-derivative}
\end{equation}
Hence $K$ decreases on $(1/\sqrt3,1]$, has its unique minimum at the regular tetrahedron $t=1$, and increases for $t>1$. Moreover,
\begin{equation}
 K(1)=6\arccos\frac13,
 \qquad
 \lim_{t\downarrow1/\sqrt3}K(t)=3\pi,
 \qquad
 \lim_{t\to\infty}K(t)=\frac{5\pi}{2}.
 \label{eq:kite-limits}
\end{equation}
The short-kite branch $t\in(1/\sqrt3,1]$ alone realizes every value in
\begin{equation}
 \{\Sigma(T):T\in\Kclass\}
 =\left[6\arccos\frac13,\,3\pi\right),
 \label{eq:kite-range}
\end{equation}
whereas the long-kite branch $t\ge1$ has range $[6\arccos(1/3),5\pi/2)$. The lower endpoint is attained at $t=1$, while the flattening short-kite branch approaches the unattained upper endpoint $3\pi$.

\begin{corollary}[Exact orthocentric range]\label{cor:orth-range}
The set of dihedral angle sums of nondegenerate orthocentric tetrahedra is
\[
 \{\Sigma(T):T\in\Oclass\}
 =\left[6\arccos\frac13,\,3\pi\right),
\]
and the lower endpoint is attained only by the regular tetrahedron.
\end{corollary}

\begin{proof}
Theorem~\ref{thm:orth-lower} gives the lower bound and its equality case, \eqref{eq:universal} gives the strict upper bound, and $\Kclass\subset\Oclass$ together with \eqref{eq:kite-range} realizes every intermediate value.
\end{proof}

\section{The isogonic lower bound and a dispersion defect}
We now turn to the main new step. Let $T=A_1A_2A_3A_4$ be isogonic. Let $O$ and $r$ denote its incenter and inradius, and let $p_i>0$ be the tangent length from $A_i$ to the insphere. Brown~\cite{Brown1924} records the classical $120^\circ$ contact property: at the point where the insphere touches any face, the three segments to the vertices of that face meet pairwise at $120^\circ$. The cosine law gives \eqref{eq:isogonicmetric}; the same contact geometry appears in Hajja--Walker~\cite{HajjaWalker2004}.

For later use put
\begin{equation}
 q_{ij}:=p_i^2+p_ip_j+p_j^2,
 \qquad
 E_p:=\sum_{1\le r<s\le4}p_rp_s,
 \qquad
 \Pi:=p_1p_2p_3p_4.
 \label{eq:qpEp}
\end{equation}
The metric representation is complete, not merely necessary. For any positive quadruple $(p_1,p_2,p_3,p_4)$, each face has the planar $120^\circ$ realization used below and hence positive area. A direct Cayley--Menger evaluation gives
\begin{equation}
 V^2=\frac{\Pi E_p}{16}>0,
 \label{eq:positive-volume}
\end{equation}
so the six lengths define a nondegenerate Euclidean tetrahedron; Hajja--Hammoudeh~\cite[Theorem~4]{HajjaHammoudeh2014} identify the resulting tetrahedron as isogonic. Thus positivity of the four parameters is the complete realization condition.

For an edge $ij$, let $\vartheta_{ij}\in(0,\pi)$ denote the interior dihedral angle at that edge, namely the angle of the wedge lying inside $T$. We first derive the required half-angle identity directly from the contact geometry.

\begin{lemma}[Geometric half-angle formula]\label{lem:geometric-half}
For complementary indices $\{i,j,k,l\}=\{1,2,3,4\}$,
\begin{equation}
 r^2=\frac{3\Pi}{4E_p},
 \label{eq:inradius-isogonic}
\end{equation}
and
\begin{equation}
 \tan^2\frac{\vartheta_{ij}}2
 =\frac{p_kp_l\,q_{ij}}{p_ip_j\,E_p}.
 \label{eq:half-p}
\end{equation}
\end{lemma}

\begin{proof}
Let $F_i$ be the contact point on the face opposite $A_i$, and write
\[
 H_i:=p_jp_k+p_jp_l+p_kp_l
 \qquad(\{i,j,k,l\}=\{1,2,3,4\}).
\]
The three $120^\circ$ angles at $F_i$ decompose the face into three triangles and give
\begin{equation}
 \Delta_i=\frac{\sqrt3}{4}H_i.
 \label{eq:face-area-contact}
\end{equation}
The barycentric weight of $A_i$ at the incenter $O$ is the ratio of the volume of the subtetrahedron with apex $O$ and base the face opposite $A_i$ to the volume of $T$. All four subtetrahedra have the same altitude $r$, so their volumes, and hence the barycentric weights, are proportional to the opposite-face areas $\Delta_i$. If $v_i:=\overrightarrow{OA_i}$, then
\[
 O=\sum_{i=1}^4\lambda_iA_i,
 \qquad
 \lambda_i=\frac{\Delta_i}{\sum_m\Delta_m},
\]
and therefore
\begin{equation}
 \sum_{i=1}^4H_iv_i=0.
 \label{eq:incenter-vector}
\end{equation}
Moreover, $p_i$ is the tangent length from $A_i$ to the insphere, so
\begin{equation}
 |v_i|^2=p_i^2+r^2.
 \label{eq:vi-length}
\end{equation}
Using $|v_i-v_j|^2=q_{ij}$, we also obtain, for $i\ne j$,
\begin{equation}
 v_i\!\cdot v_j=r^2-\frac12p_ip_j.
 \label{eq:vi-dot}
\end{equation}

Let
\[
 e_3:=\sum_{1\le i<j<k\le4}p_ip_jp_k.
\]
The elementary identities
\[
 \sum_iH_i=2E_p,
 \qquad
 \sum_ip_iH_i=3e_3,
 \qquad
 e_3-p_jH_j=\frac{\Pi}{p_j}
\]
will be used below. Taking the scalar product of \eqref{eq:incenter-vector} with $v_j$ and separating the term $i=j$ gives
\begin{align}
 0
 &=H_j|v_j|^2+\sum_{i\ne j}H_i(v_i\!\cdot v_j)\notag\\
 &=H_j(p_j^2+r^2)
   +\sum_{i\ne j}H_i\left(r^2-\frac12p_ip_j\right)\notag\\
 &=r^2\sum_iH_i+H_jp_j^2
   -\frac{p_j}{2}\left(\sum_ip_iH_i-p_jH_j\right)\notag\\
 &=2E_pr^2+H_jp_j^2
   -\frac{p_j}{2}\bigl(3e_3-p_jH_j\bigr).             \label{eq:inradius-expand}
\end{align}
Since $e_3=p_jH_j+\Pi/p_j$,
\[
 3e_3-p_jH_j=2p_jH_j+\frac{3\Pi}{p_j}.
\]
Substitution makes the two $H_jp_j^2$ terms cancel and leaves
\[
 0=2E_pr^2-\frac32\Pi,
\]
which proves \eqref{eq:inradius-isogonic}. This also recovers the inradius formula obtained in Hajja--Walker~\cite{HajjaWalker2004}.

Now fix the edge $A_iA_j$, with complementary indices $\{i,j,k,l\}=\{1,2,3,4\}$, and choose the contact point $F_k$ on the adjacent face opposite $A_k$. The triangle $A_iF_kA_j$ has sides $p_i,p_j$ enclosing an angle $120^\circ$, so $A_iA_j=\sqrt{q_{ij}}$ and the perpendicular distance $s_{ij}$ from $F_k$ to the edge is
\begin{equation}
 s_{ij}=\frac{\sqrt3\,p_ip_j}{2\sqrt{q_{ij}}}.
 \label{eq:sij}
\end{equation}
The other adjacent contact point $F_l$ gives the same distance, because $A_iF_lA_j$ has the same two side lengths and the same included angle.

Because $OF_k$ is normal to the face containing $A_iA_j$, we have $OF_k\perp A_iA_j$. Hence $F_k$ lies in the plane through $O$ perpendicular to $A_iA_j$; see Figure~\ref{fig:kite-crosssection}(b). The intersection of the tetrahedron with this plane is the interior wedge between the two adjacent faces, so its angle is precisely the interior dihedral angle $\vartheta_{ij}$, not its supplement. The insphere becomes a circle of radius $r$ tangent to both boundary lines. Since $O$ is equidistant from the two faces, it lies on the internal bisector of the wedge, and the right triangle $POF_k$ gives
\begin{equation}
 s_{ij}=r\cot\frac{\vartheta_{ij}}2.
 \label{eq:wedge-cross-section}
\end{equation}
Combining \eqref{eq:inradius-isogonic}, \eqref{eq:sij}, and \eqref{eq:wedge-cross-section} yields \eqref{eq:half-p}.
\end{proof}

\begin{remark}[Independent normal-vector check]
A direct cofactor calculation for the edge Gram matrix gives
\[
 \cos\vartheta_{ij}
 =1-\frac{2p_kp_lq_{ij}}{H_kH_l}.
\]
The identity
\[
 H_kH_l-p_kp_lq_{ij}=p_ip_jE_p
\]
then converts $\tan^2(\vartheta_{ij}/2)=(1-\cos\vartheta_{ij})/(1+\cos\vartheta_{ij})$ into \eqref{eq:half-p}. This provides an algebraic check of the interior-angle convention used in the cross-sectional proof.
\end{remark}

The reciprocal form of Lemma~\ref{lem:geometric-half} is especially simple.

\begin{lemma}[Reciprocal half-angle form]\label{lem:reciprocal}
Set
\begin{equation}
 x_i:=\frac1{p_i},
 \qquad
 E:=\sum_{1\le r<s\le4}x_rx_s,
 \qquad
 Q_{ij}:=x_i^2+x_ix_j+x_j^2.
 \label{eq:xEQ}
\end{equation}
Then for every edge $ij$,
\begin{equation}
 \tan^2\frac{\vartheta_{ij}}2=\frac{Q_{ij}}{E}.
 \label{eq:reciprocal-half}
\end{equation}
\end{lemma}

\begin{proof}
For complementary $\{i,j,k,l\}=\{1,2,3,4\}$, substitute $p_r=1/x_r$ into \eqref{eq:half-p}. Since
\[
 E_p=\frac{E}{x_1x_2x_3x_4},
 \qquad
 q_{ij}=\frac{Q_{ij}}{x_i^2x_j^2},
\]
all remaining factors cancel.
\end{proof}

\begin{remark}[Classical opposite-edge characterizations]\label{rem:angle-characterizations}
In the Russian literature on tetrahedral geometry, Ponarin~\cite{Ponarin2009} records the characteristic isogonic relation
\begin{equation}
 \cos\frac{\vartheta_{12}}2\cos\frac{\vartheta_{34}}2
 =\cos\frac{\vartheta_{13}}2\cos\frac{\vartheta_{24}}2
 =\cos\frac{\vartheta_{14}}2\cos\frac{\vartheta_{23}}2.
 \label{eq:ponarin-isogonic}
\end{equation}
For comparison, Sharygin~\cite[Problem~312(e)]{Sharygin1986} lists
\begin{equation}
 \cos\vartheta_{12}\cos\vartheta_{34}
 =\cos\vartheta_{13}\cos\vartheta_{24}
 =\cos\vartheta_{14}\cos\vartheta_{23}
 \label{eq:sharygin-orthocentric}
\end{equation}
among the equivalent criteria for orthocentricity. Thus the same three opposite-edge pairings occur naturally in both classes: the isogonic characterization uses half-angles, whereas the orthocentric characterization uses full dihedral angles. These identities do not imply the common angle-sum bound, but they explain why opposite-edge pairing is natural in the isogonic proof below.
\end{remark}

Group the six angles into the three opposite-edge pairs
\[
 12\mid34,
 \qquad
 13\mid24,
 \qquad
 14\mid23.
\]
For one such partition $ij\mid kl$, put
\begin{equation}
 \Phi_{ij\mid kl}:=\frac{\vartheta_{ij}+\vartheta_{kl}}2.
 \label{eq:Phi}
\end{equation}
Since
\[
 Q_{ij}=(x_i+x_j)^2-x_ix_j<(x_i+x_j)^2,
 \qquad
 Q_{kl}<(x_k+x_l)^2,
\]
and
\[
 E=(x_i+x_j)(x_k+x_l)+x_ix_j+x_kx_l,
\]
we have
\begin{equation}
 \sqrt{Q_{ij}Q_{kl}}<(x_i+x_j)(x_k+x_l)<E.
 \label{eq:product-less-E}
\end{equation}
By Lemma~\ref{lem:reciprocal},
\[
 \tan\frac{\vartheta_{ij}}2\tan\frac{\vartheta_{kl}}2
 =\frac{\sqrt{Q_{ij}Q_{kl}}}{E}<1.
\]
Both half-angles lie in $(0,\pi/2)$, so $0<\Phi_{ij\mid kl}<\pi/2$. Define
\begin{equation}
 \psi_{ij\mid kl}:=\frac\pi2-\Phi_{ij\mid kl}\in(0,\pi/2).
 \label{eq:psi}
\end{equation}
For brevity set
\[
 \psi_1:=\psi_{12\mid34},
 \qquad
 \psi_2:=\psi_{13\mid24},
 \qquad
 \psi_3:=\psi_{14\mid23}.
\]
The tangent addition formula gives
\begin{equation}
 \tan\psi_{ij\mid kl}
 =\frac{E-\sqrt{Q_{ij}Q_{kl}}}
 {\sqrt E\,(\sqrt{Q_{ij}}+\sqrt{Q_{kl}})}.
 \label{eq:tanpsi-exact}
\end{equation}

The proof first bounds each complementary opposite-edge angle. Summation removes the dependence on the partition, while an exact dispersion identity for the reciprocal tangent parameters produces the regular-tetrahedron constant and a quantitative defect term.

\begin{theorem}[Isogonic lower bound with dispersion defect]\label{thm:isogonic-lower}
Let
\begin{equation}
 S:=x_1+x_2+x_3+x_4,
 \qquad
 D:=\sum_{1\le i<j\le4}(x_i-x_j)^2,
 \qquad
 \delta:=\frac{D}{3S^2}.
 \label{eq:dispersion}
\end{equation}
Then $0\le\delta<1$ and
\begin{equation}
 \Sigma(T)\ge
 3\pi-6\arctan\!\left(
 \frac{1}{2\sqrt2}\sqrt{1-\delta}
 \right).
 \label{eq:defect-bound}
\end{equation}
Equality in \eqref{eq:defect-bound} holds if and only if $T$ is regular. In particular,
\begin{equation}
 \Sigma(T)\ge6\arccos\frac13,
 \label{eq:isogonic-lower}
\end{equation}
with equality if and only if $T$ is regular.
\end{theorem}

\begin{proof}
The exact local defect identity
\begin{equation}
 x^2+xy+y^2-\frac34(x+y)^2=\frac14(x-y)^2
 \label{eq:local-defect}
\end{equation}
gives
\begin{equation}
 \sqrt{x^2+xy+y^2}\ge\frac{\sqrt3}{2}(x+y),
 \qquad x,y>0,
 \label{eq:Q-lower}
\end{equation}
with equality exactly when $x=y$. For a fixed opposite-edge partition $ij\mid kl$, put
\[
 u_{ij\mid kl}:=x_ix_j+x_kx_l.
\]
Since $(x_i+x_j)(x_k+x_l)=E-u_{ij\mid kl}$, inequality \eqref{eq:Q-lower} gives
\[
 \sqrt{Q_{ij}Q_{kl}}\ge\frac34(E-u_{ij\mid kl}),
 \qquad
 \sqrt{Q_{ij}}+\sqrt{Q_{kl}}\ge\frac{\sqrt3}{2}S.
\]
Therefore
\begin{equation}
 E-\sqrt{Q_{ij}Q_{kl}}
 \le\frac14(E+3u_{ij\mid kl}),
 \label{eq:numerator-bound}
\end{equation}
and
\begin{equation}
 \sqrt E\,(\sqrt{Q_{ij}}+\sqrt{Q_{kl}})
 \ge\frac{\sqrt3}{2}S\sqrt E.
 \label{eq:denominator-bound}
\end{equation}
Substitution into \eqref{eq:tanpsi-exact} yields
\begin{equation}
 \tan\psi_{ij\mid kl}
 \le\frac{E+3u_{ij\mid kl}}{2\sqrt3\,S\sqrt E}.
 \label{eq:tanpsi-bound}
\end{equation}
Because the numerator is bounded from above and the positive denominator from below, equality in \eqref{eq:tanpsi-bound} requires equality in both applications of \eqref{eq:Q-lower}; it therefore holds for a fixed partition exactly when $x_i=x_j$ and $x_k=x_l$.

Since
\[
 u_{12\mid34}+u_{13\mid24}+u_{14\mid23}=E,
\]
summing \eqref{eq:tanpsi-bound} gives
\begin{equation}
 \sum_{m=1}^3\tan\psi_m
 \le\frac{3E+3E}{2\sqrt3\,S\sqrt E}
 =\frac{\sqrt{3E}}{S}.
 \label{eq:sumtan-first}
\end{equation}
Equality here requires equality for all three partitions and therefore forces $x_1=x_2=x_3=x_4$.

The global symmetry defect is exact:
\begin{equation}
 D=4\sum_{i=1}^4x_i^2-S^2,
 \qquad
 E=\frac{S^2-\sum_i x_i^2}{2}
   =\frac{3S^2-D}{8}.
 \label{eq:global-defect}
\end{equation}
Since $E>0$, this implies $0\le\delta<1$. Combining \eqref{eq:sumtan-first} and \eqref{eq:global-defect},
\begin{equation}
 \frac13\sum_{m=1}^3\tan\psi_m
 \le\frac{1}{2\sqrt2}\sqrt{1-\delta}.
 \label{eq:average-tan-defect}
\end{equation}
The function $\tan x$ is strictly convex on $(0,\pi/2)$. Hence Jensen's inequality and \eqref{eq:average-tan-defect} give
\[
 \tan\!\left(\frac{\psi_1+\psi_2+\psi_3}{3}\right)
 \le\frac13\sum_{m=1}^3\tan\psi_m
 \le\frac{1}{2\sqrt2}\sqrt{1-\delta}.
\]
Since $\tan$ is increasing,
\begin{equation}
 \psi_1+\psi_2+\psi_3
 \le3\arctan\!\left(
 \frac{1}{2\sqrt2}\sqrt{1-\delta}
 \right).
 \label{eq:psi-defect}
\end{equation}
Finally,
\[
 \frac{\Sigma(T)}2
 =\Phi_1+\Phi_2+\Phi_3
 =\frac{3\pi}{2}-(\psi_1+\psi_2+\psi_3),
\]
which proves \eqref{eq:defect-bound}.

If equality holds in \eqref{eq:defect-bound}, then equality must hold in both Jensen's inequality and \eqref{eq:sumtan-first}. The latter already forces $x_1=x_2=x_3=x_4$, so $\delta=0$ and $T$ is regular. Conversely, for the regular tetrahedron every inequality above is an equality. The sharp lower bound \eqref{eq:isogonic-lower} follows from $\delta\ge0$ and
\[
 3\pi-6\arctan\frac{1}{2\sqrt2}
 =3\pi-6\arcsin\frac13
 =6\arccos\frac13.
\]
\end{proof}

\begin{remark}[Near equality in reciprocal-parameter space]
Estimate \eqref{eq:defect-bound} is a parameter-space defect inequality, not by itself a stability theorem in a chosen geometric metric on similarity classes. It nevertheless gives a quantitative reverse implication. Put $a:=\arcsin(1/3)$. If $0\le\varepsilon<6a$ and
\[
 \Sigma(T)\le\Sigma_{\rm reg}+\varepsilon,
\]
then
\begin{equation}
 \delta\le1-8\tan^2\!\left(a-\frac{\varepsilon}{6}\right).
 \label{eq:near-equality-delta}
\end{equation}
After the similarity normalization $S=1$, identity \eqref{eq:global-defect} gives
\[
 \sum_{i=1}^4\left(x_i-\frac14\right)^2
 =\frac D4
 =\frac{3\delta}{4}.
\]
Thus near equality forces the normalized reciprocal tangent parameters to be close to equality.
\end{remark}

\begin{proof}[Proof of Theorem~\ref{thm:isogonic-main}]
The lower bound and equality case are contained in Theorem~\ref{thm:isogonic-lower}. The universal bound \eqref{eq:universal} gives $\Sigma(T)<3\pi$. Finally, Proposition~\ref{prop:intersection} and the short-kite branch in \eqref{eq:kite-range} realize every value in
\[
 \left[6\arccos\frac13,\,3\pi\right).
\]
\end{proof}

\begin{remark}[Relation to Hammoudeh's ordered reduction]\label{rem:ordered-reduction}
Hammoudeh~\cite[Appendix~A, Conjecture~A.2]{Hammoudeh2026} proposes the stronger inequality obtained by ordering $0<p_1\le p_2\le p_3\le p_4$, setting $m=(p_1+p_2+p_3)/3$, and comparing
\[
 \Sigma_{\rm igo}(p_1,p_2,p_3,p_4)
 \quad\text{with}\quad
 \Sigma_{\rm igo}(m,m,m,p_4).
\]
The conjecture asserts that the second quantity is no larger, with equality only when $p_1=p_2=p_3$. Theorem~\ref{thm:isogonic-lower} proves the sharp lower bound and the dispersion defect estimate without establishing this stronger monotone comparison with a specified kite.
\end{remark}

\section{Conclusion and an open question}
Corollary~\ref{cor:orth-range} and Theorem~\ref{thm:isogonic-main} give
\[
 \{\Sigma(T):T\in\Oclass\}
 =\{\Sigma(T):T\in\Iclass\}
 =\{\Sigma(T):T\in\Kclass\}
 =\left[6\arccos\frac13,\,3\pi\right).
\]
The role of the common kite family is precise: it supplies the regular lower-endpoint example and, through the short-kite branch, realizes every value up to the unattained endpoint $3\pi$. It does not prove either lower bound; those arise from two different mechanisms.

For orthocentric tetrahedra the decisive input is the four-corner conservation law \eqref{eq:polar-sum} combined with a supporting line. For isogonic tetrahedra the reciprocal tangent parameters carry the scale-invariant defect $\delta$ in \eqref{eq:dispersion}, and Theorem~\ref{thm:isogonic-lower} converts that defect into an explicit angle-sum gap. The near-equality estimate \eqref{eq:near-equality-delta} controls reciprocal-parameter dispersion, but it does not construct a deformation in the geometric space of tetrahedra.

This distinction motivates a natural open question: does there exist a geometrically natural deformation from an arbitrary isogonic tetrahedron to a kite along which the dihedral angle sum is nonincreasing? The ordered reduction in Remark~\ref{rem:ordered-reduction} asks for a specific comparison of this kind. Such a result would complement the present defect inequality by supplying a monotone geometric route to the extremal kite family.


\begin{thebibliography}{99}
\small

\bibitem{BrandtsCihangirKrizek2015}
J.~Brandts, A.~Cihangir, and M.~K\v{r}\'i\v{z}ek,
\emph{Tight bounds on angle sums of nonobtuse simplices},
Appl. Math. Comput. \textbf{267} (2015), 397--408.
\href{https://doi.org/10.1016/j.amc.2015.02.035}{doi:10.1016/j.amc.2015.02.035}.

\bibitem{Brown1924}
B.~H. Brown,
\emph{A theorem on isogonal tetrahedra},
Amer. Math. Monthly \textbf{31} (1924), no.~8, 371--375.

\bibitem{EdmondsHajjaMartini2005}
A.~L. Edmonds, M.~Hajja, and H.~Martini,
\emph{Orthocentric simplices and their centers},
Results Math. \textbf{47} (2005), 266--295.
\href{https://doi.org/10.1007/BF03323029}{doi:10.1007/BF03323029}.

\bibitem{Eriksson1990}
F.~Eriksson,
\emph{On the measure of solid angles},
Math. Mag. \textbf{63} (1990), no.~3, 184--187.
\href{https://doi.org/10.1080/0025570X.1990.11977515}{doi:10.1080/0025570X.1990.11977515}.

\bibitem{Gaddum1952}
J.~W. Gaddum,
\emph{The sums of the dihedral and trihedral angles in a tetrahedron},
Amer. Math. Monthly \textbf{59} (1952), 370--371.
\href{https://doi.org/10.1080/00029890.1952.11988143}{doi:10.1080/00029890.1952.11988143}.

\bibitem{HajjaHammoudeh2014}
M.~Hajja and I.~Hammoudeh,
\emph{The sum of measures of the angles of a simplex},
Beitr. Algebra Geom. \textbf{55} (2014), 453--470.
\href{https://doi.org/10.1007/s13366-013-0160-8}{doi:10.1007/s13366-013-0160-8}.

\bibitem{HajjaHammoudehHayajneh2015}
M.~Hajja, I.~Hammoudeh, and M.~Hayajneh,
\emph{Kites as the only doubly special simplices},
Beitr. Algebra Geom. \textbf{56} (2015), 269--277.
\href{https://doi.org/10.1007/s13366-014-0204-8}{doi:10.1007/s13366-014-0204-8}.

\bibitem{HajjaHayajneh2014}
M.~Hajja and M.~Hayajneh,
\emph{Impurity of the corner angles in certain special families of simplices},
J. Geom. \textbf{105} (2014), 539--560.
\href{https://doi.org/10.1007/s00022-014-0219-1}{doi:10.1007/s00022-014-0219-1}.

\bibitem{HajjaHayajnehMartini2016}
M.~Hajja, M.~Hayajneh, and H.~Martini,
\emph{More characterizations of certain special families of simplices},
Results Math. \textbf{69} (2016), 23--47.
\href{https://doi.org/10.1007/s00025-015-0456-0}{doi:10.1007/s00025-015-0456-0}.

\bibitem{HajjaWalker2004}
M.~Hajja and P.~Walker,
\emph{The inspherical Gergonne center of a tetrahedron},
J. Geom. Graph. \textbf{8} (2004), no.~1, 23--32.

\bibitem{Hammoudeh2026}
I.~Hammoudeh,
\emph{Dihedral-Angle Sums in Orthocentric, Power-$k$, and Isodynamic Tetrahedra},
Version~2.0.0, 6 September 2026; \href{https://ismaa3iil.fyi/orthocentric-powerk-dihedral-sums/}{author's website} (accessed 19 October 2026).

\bibitem{Jarden1941}
D.~Jarden (D.~Juzuk),
\emph{On the sum of the dihedral and trihedral angles in a tetrahedron} (in Hebrew),
Technika Umada (Technics and Science) \textbf{4} (1941), no.~5 (26), 8--9.

\bibitem{Katsuura2020}
H.~Katsuura,
\emph{Solid angle sum of a tetrahedron},
J. Geom. Graph. \textbf{24} (2020), no.~1, 29--34.

\bibitem{LiuJoe1994}
A.~Liu and B.~Joe,
\emph{Relationship between tetrahedron shape measures},
BIT \textbf{34} (1994), no.~2, 268--287.
\href{https://doi.org/10.1007/BF01955874}{doi:10.1007/BF01955874}.

\bibitem{KorotovKrizek2026}
S.~Korotov and M.~K\v{r}\'i\v{z}ek,
\emph{Two-sided bounds for dihedral angle sums of path and 4-ball tetrahedra},
arXiv:2601.08304v1 [math.MG], 2026.

\bibitem{Mandan1961}
S.~R. Mandan,
\emph{Isodynamic and isogonic simplexes},
Ann. Mat. Pura Appl. (4) \textbf{53} (1961), 45--55.
\href{https://doi.org/10.1007/BF02417784}{doi:10.1007/BF02417784}.

\bibitem{Ponarin2009}
Ya.~P. Ponarin,
\emph{Elementarnaya geometriya. Tom 3: Treugolniki i tetraedry [Elementary Geometry, Vol.~3: Triangles and Tetrahedra]},
MCCME, Moscow, 2009 (in Russian).

\bibitem{Sharygin1986}
I.~F. Sharygin,
\emph{Problems in Solid Geometry},
Mir Publishers, Moscow, 1986; revised from the 1984 Russian edition.

\bibitem{SkopetsPonarin1974}
Z.~A. Skopets and Ya.~P. Ponarin,
\emph{Geometriya tetraedra i ego elementov [Geometry of the Tetrahedron and Its Elements]},
Yaroslavl, 1974 (in Russian).

\end{thebibliography}
\end{document}